# DEDUCIBILITY THEOREMS IN BOOLEAN LOGIC

Florentin Smarandache
University of New Mexico
200 College Road
Gallup, NM 87301, USA
E-mail: smarand@unm.edu

**ABSTRACT**
In this paper we give two theorems from the Propositional Calculus of the Boolean Logic with their consequences and applications and we prove them axiomatically.

## §1. THEOREMS, CONSEQUENCES

In the beginning I shall put forward the axioms of the Propositional Calculus.

I.  a) $\vdash A \supset (B \supset A)$,
    b) $\vdash (A \supset (B \supset C)) \supset ((A \supset B) \supset (A \supset C))$.

II. a) $\vdash A \wedge B \supset A$,
    b) $\vdash A \wedge B \supset B$,
    c) $\vdash (A \supset B) \supset ((A \supset C) \supset (A \supset B \wedge C))$.

III. a) $\vdash A \supset A \vee B$,
     b) $\vdash B \supset A \vee B$,
     c) $\vdash (A \supset C) \supset ((B \supset C) \supset (A \vee B \supset C))$.

IV. a) $\vdash (A \supset B) \supset (\overline{B} \supset \overline{A})$,
    b) $\vdash A \supset \overline{\overline{A}}$,
    c) $\vdash \overline{\overline{A}} \supset A$.

**THEOREMS.** If: $\vdash A_i \supset B_i, i = \overline{1,n}$, then

1) $\vdash A_1 \wedge A_2 \wedge ... \wedge A_n \supset B_1 \wedge B_2 \wedge ... \wedge B_n$,
2) $\vdash A_1 \vee A_2 \vee ... \vee A_n \supset B_1 \vee B_2 \vee ... \vee B_n$.

*Proof:*

It is made by complete induction. For $n = 1$: $\vdash A_1 \supset B_1$, which is true from the given hypothesis. For $n = 2$: hypotheses $\vdash A_1 \supset B_1$, $\vdash A_2 \supset B_2$; let's show that $\vdash A_1 \wedge A_2 \supset B_1 \wedge B_2$. We use the axiom II, c) replacing $A \to A_1 \wedge A_2$, $B \to B_1$, $C \to B_2$, it results:

(1)  $\vdash (A_1 \wedge A_2 \supset B_1) \supset ((A_1 \wedge A_2 \supset B_2) \supset (A_1 \wedge A_2 \supset B_1 \wedge B_2))$.

We use the axiom II, a) replacing $A \to A_1$, $B \to A_2$; we have $\vdash A_1 \wedge A_2 \supset A_1$. But $\vdash A_1 \supset B_1$ (hypothesis) applying the syllogism rule, it results $\vdash A_1 \wedge A_2 \supset B_1$. Analogously, using the axiom II, b), we have $\vdash A_1 \wedge A_2 \supset B_2$. We know that $\vdash A_1 \wedge A_2 \supset B_i$, $i = 1, 2$, are deducible, then applying in (I) inference rule twice, we have $\vdash A_1 \wedge A_2 \supset B_1 \wedge B_2$.

We suppose it's true for $n$; let's prove that for $n+1$ it is true. In $\vdash A_1 \wedge A_2 \supset B_1 \wedge B_2$ replacing $A_1 \to A_1 \wedge ... \wedge A_n$, $A_2 \to A_{n+1}$, $B_1 \to B_1 \wedge ... \wedge B_n$, $B_2 \to B_{n+1}$ and using induction hypothesis it results $\vdash A_1 \wedge ... \wedge A_n \wedge A_{n+1} \supset B_1 \wedge ... \wedge B_n \wedge B_{n+1}$ and item 1) from the Theorem is proved.

2) It is made by induction. For $n=1$; if $\vdash A_1 \supset B_1$, then of course $\vdash A_1 \supset B_1$. For $n=2$: if $\vdash A_1 \supset B_1$ and $\vdash A_2 \supset B_2$, then $\vdash A_1 \vee A_2 \supset B_1 \vee B_2$.

We use axiom III, c) replacing $A \to A_1$, $B \to A_2$, $C \to B_1 \vee B_2$ we get

(2) $\quad \vdash (A_1 \supset B_1 \vee B_2) \supset ((A_2 \supset B_1 \vee B_2) \supset (A_1 \vee A_2 \supset B_1 \vee B_2))$.

Let's show that $\vdash A_1 \supset B_1 \vee B_2$. We use the axiom III, a) replacing $A \to B_1$, $B \to B_2$ we get $\vdash B_1 \supset B_1 \vee B_2$ and we know from the hypothesis $A_1 \quad B_1$. Applying the syllogism we get $\vdash A_1 \supset B_1 \vee B_2$.

In the axiom III, b) replacing $A \to B_1$, $B \to B_2$, we get $\vdash B_2 \supset B_1 \vee B_2$. But $\vdash A_2 \supset B_2$ (from the hypothesis), applying the syllogism we get $\vdash A_2 \supset B_1 \vee B_2$. Applying the inference rule twice in (2) we get $\vdash A_1 \vee A_2 \supset B_1 \vee B_2$.

Suppose it's true for $n$ and let's show that for $n+1$ it is true. Replace in $\vdash A_1 \vee A_2 \supset B_1 \vee B_2$ (true formula if $\vdash A_1 \supset B_1$ and $\vdash A_2 \supset B_2$) $A_1 \to A_1 \vee ... \vee A_n$, $A_2 \to A_{n+1}$, $B_1 \to B_1 \vee ... \vee B_n$, $B_2 \to B_{n+1}$. From induction hypothesis it results $\vdash A_1 \vee ... \vee A_n \vee A_{n+1} \supset B_1 \vee ... \vee B_n \vee B_{n+1}$ and the theorem is proved.

**CONSEQUENCES.**

1°) If $\vdash A_i \supset B$, $i = \overline{1,n}$ then $\vdash A_1 \wedge ... \wedge A_n \supset B$.

2°) If $\vdash A_i \supset B$, $i = \overline{1,n}$, then $\vdash A_1 \vee ... \vee A_n \supset B$.

*Proof:* 1°) Using 1) from the theorem, we get

(3) $\quad \vdash A_1 \wedge ... \wedge A_n \supset B \wedge ... \wedge B$ ($n$ times).

In axiom II, a) we replace $A \to B$, $B \to B \wedge ... \wedge B$ ($n-1$ times), and we get

(4) $\quad \vdash B \wedge ... \wedge B \supset B$ (n times).

From (3) and (4) by means of the syllogism rule we get $\vdash A_1 \wedge ... \wedge A_n \supset B$.

2°) Using 2) from theorem, we get $\vdash A_1 \vee ... \vee A_n \supset B \vee ... \vee B$ ($n$ times).

**LEMMA.** $\vdash B \vee ... \vee B \supset B$ ($n$ times), $n \geq 1$.

*Proof:*

It is made by induction. For $n=1$, obvious. For $n=2$: in axiom III, c) we replace $A \to B$, $C \to B$ and we get $\vdash (B \supset B) \supset ((B \supset B) \supset (B \vee B \supset B))$. Applying the inference rule twice we get $\vdash B \vee B \supset B$.

Suppose for $n$ that the formula is deducible, let's prove that is for $n+1$.

We proved that $\vdash B \supset B$. In axiom III, c) we replace $A \to B \vee ... \vee B$ ($n$ times), $C \to B$, and we get $\vdash (B \vee ... \vee B \supset B) \supset ((B \supset B) \supset (B \vee ... \vee B \supset B))$ ($n$ times). Applying two times the interference rule, we get $\vdash B \vee ... \vee B \supset B$ ($n+1$ times) so lemma is proved.

From $\vdash A_1 \vee ... \vee A_n \supset B \vee ... \vee B$ ($n$ times) and applying the syllogism rule, from lemma we get $\vdash A_1 \vee ... \vee A_n \supset B$.

3°) $\vdash A \wedge ... \wedge A \supset A$ ($n$ times)
4°) $\vdash A \vee ... \vee A \supset A$ ($n$ times).

Previously we proved, replacing in Consequence 1°) and 2°), $B \to A$. Analogously, the consequences are proven:

5°) If $\vdash A \supset B_i, i = \overline{1,n}$, then $\vdash A \supset B_1 \wedge ... \wedge B_n$.

6°) If $\vdash A \supset B_i, i = \overline{1,n}$, then $\vdash A \supset B_1 \vee ... \vee B_n$.

Analogously,

7°) $\vdash A \supset A \wedge ... \wedge A$ ($n$ times)
8°) $\vdash A \supset A \vee ... \vee A$ ($n$ times)
9°) $\vdash A_1 \wedge ... \wedge A_n \supset A_1 \vee ... \vee A_n$.

*Proof:*

Method I. It is initially proved by induction: $\vdash A_1 \wedge ... \wedge A_n \supset A_i$, $i = \overline{1,n}$ and 2) is applied from the Theorem.

Method II. It is proven by induction that: $\vdash A_i \supset A_1 \wedge ... \wedge A_n$, $i = \overline{1,n}$ and then 1) is applied from the Theorem.

10°) If $\vdash A_i \supset B_i$, $i = \overline{1,n}$, then $\vdash A_1 \wedge ... \wedge A_n \supset B_1 \vee ... \vee B_n$.

*Proof:*

Method I. Using 1) from the Theorem, it results:

(5)    $\vdash A_1 \wedge ... \wedge A_n \supset B_1 \wedge ... \wedge B_n$

We apply the Consequence 9°) where we replace $A_i \to B_i$, $i = \overline{1,n}$ and results:

(6)    $\vdash B_1 \wedge ... \wedge B_n \supset B_1 \vee ... \vee B_n$.

From (5) and (6), applying the syllogism rule we get 10°).

Method II. We firstly use the Consequence 9°) and then 2) from the Theorem and so we obtain the Consequence 10°).

## §2. APPLICATIONS AND REMARKS ON THEOREMS

The theorems are used in order to prove the formulae of the shape:

$\vdash A_1 \wedge ... \wedge A_p \supset B_1 \wedge ... \wedge B_r$

$\vdash A_1 \vee ... \vee A_p \supset B_1 \vee ... \vee B_r$, where $p, r \in \mathbb{N}^*$

It is proven that $\vdash A_i \supset B_j$, i.e.

$\forall i \in \overline{1,p}$, $\exists j_0 \in \overline{1,r}$, $j_0 = j_0(i)$, $\vdash A_i \supset B_{j_0}$

and

$\forall j \in \overline{1,r}$, $\exists i_0 \in \overline{1,p}$, $i_0 = i_0(j)$, $\vdash A_{i_0} \supset B_j$.

**EXAMPLES**: The following formulas are deducible:

(i)     $\vdash A \supset (A \vee B) \wedge (B \supset A)$,

(ii)    $\vdash (A \wedge B) \vee C \supset A \vee B \vee C$,

(iii)   $\vdash A \wedge C \supset A \vee C$.

***Solution:***

(i) We have $\vdash A \supset A \vee B$ and $\vdash A \supset (B \supset A)$ (axiom III, a) and I, a)) and according to 1) from Theorem it results (i).

(ii) From $\vdash A \supset (B \supset A)$, $\vdash A \wedge B \supset B$, $\vdash C \supset C$ and Theorem 1), we have (ii).

(iii) Method I. From $\vdash A \wedge C \supset A$, $\vdash A \wedge C \supset C$ and Theorem 2).

Method II. From $\vdash A \supset A \vee C$, $\vdash C \supset A \vee C$ and using Theorem 1).

**REMARKS.** 1) The reciprocals of Theorem 1) and 2) are not always true.

a) Counter-example for Theorem 1). The formula $\vdash A \wedge B \supset A \wedge A$ is deducible from axiom II, a), $\vdash A \wedge A \supset A$ (Consequence 7°) and the syllogism rule. But $\vdash A \supset A$ for all A, that the formula $B \supset A$ is not deducible, so the reciprocal of the Theorem 1) is false.

Counter-example for Theorem 2). The formula $\vdash A \vee A \supset A \vee B$ is deducible from Lemma, axiom III, a) and applying the syllogism rule. But $\vdash A \supset A$ for all A, that the formula $A \supset B$ is not deducible, so the reciprocal of Theorem 2) is false.

2) The reciprocals of Theorem 1) and 2) are not always true.

Counter-examples:

a) for Theorem 1): $\vdash A \supset A$ and $B \not\supset A$ results that $\vdash A \wedge B \supset A \wedge A$ so the reciprocal of Theorem 1) is false.

b) for Theorem 2): $\vdash A \supset A$ and $A \not\supset B$ results that $\vdash A \vee A \supset A \vee B$ so the reciprocal of Theorem 2) is false.

UNIVERSITATEA DIN CRAIOVA
*Facultatea de Ştiinte Exacte*
24.10.1979